\documentstyle[12pt,amsfonts]{article}

\input cyracc.def

\begin{document}
\title{The Primary Pretenders}
\author{John H.~Conway, Richard K.~Guy, W.~A.~Schneeberger \& N.~J.~A.~Sloane}
\maketitle
\setlength{\baselineskip}{1.5\baselineskip}
\thispagestyle{empty}

Perhaps the most famous theorem in number theory is Fermat's theorem.
Not Fermat's Last Theorem, of course, because that's now old hat, but
Fermat's Little Theorem:

\smallskip

\centerline{If $p$ is a prime, and $b$ is a positive integer prime
 to $p$, then $b^{p-1}\equiv 1$ (mod $p$),}

\smallskip

\noindent
which we prefer to write in the simpler form
$$b^p\equiv b\pmod p.$$

If the converse of the theorem were true, then number theory
would be a lot simpler than it is, but fortunately that is not the
case. Counterexamples to the converse of the first (and, very occasionally, the second) form of
Fermat's theorem are called {\bf pseudoprimes}. A well-known example
is $341=11\times31$, which is a pseudoprime to base 2:
$$2^{340}\equiv1\pmod{341}$$
The literature on pseudoprimes is extensive; for an introduction see section
{\bf A12} of the second author's {\it Unsolved Problems in Number Theory},
2nd edition, Springer, 1994.  D.H. Lehmer found the even pseudoprime
$161038 = 2 \cdot 73 \cdot 1103$ and N.G.W.H. Beeger showed that there were
infinitely many.

The {\bf Carmichael numbers}, such as $561=3\times11\times17$,
are counterexamples to the second form of Fermat's theorem to {\em any} base:
$$b^{561}\equiv b\pmod{561}~,~~~b=1,2, \ldots ~.$$
The second form of the theorem admits a much wider class of counterexamples
than the first, and to distinguish them from the pseudoprimes we will
call any composite number $q$ such that
$b^q\equiv b$ (mod $q$) a {\bf prime pretender} to base $b$.

We investigate $q_b$, the least prime pretender, or {\bf primary
pretender}, for the base $b$.

We will see that there are only 132 distinct primary pretenders,
and that $q_b$ is a periodic function of $b$ whose period is
 the 122-digit number

19\,\,5685843334\,\,6007258724\,\,5340037736\,\,2789820172\,\,1382933760\,\,4336734362-\\
\verb+ +\qquad\qquad
 2947386477\,\,7739548319\,\,6097971852\,\,9992599213\,\,2923650684\,\,2360439300

What is this number?  Well, it's $p!_{59}p!_9$, where $p!_k$ is
the product of the first $k$
 primes, $p_1p_2\cdots p_k$. And where do $p_{59}$ and $p_9$ come
 from? $p_{59}=277$ is the largest possible prime factor, and $p_9=23$
is the largest possible repeated prime factor, of a composite number
 less than the Carmichael number 561.

For what bases is 4 a prime pretender? If $b\equiv0$, 1, 2, 3 (mod 4),
 then $b^4\equiv0$, 1, 0, 1, so 4 is a prime pretender just for
 $b\equiv0$, 1 (mod 4).

The similar calculations mod 6 and 8 show that 6 is a prime pretender
for bases $\equiv0$ or 1 (mod 3) and that 8 is a prime pretender
 for bases $\equiv0$ or 1 (mod 8).  It follows that every number for
 which 8 is a prime pretender also has 4 as a prime pretender, so that
 8 can never be the {\it primary\/} pretender.  The calculations mod 9
 show that 9 is a prime pretender for bases  $\equiv0$, 1 or 8 (mod 9),
 which may also be described as the square roots of 0 or 1 (mod 9). 

These results can be recorded by saying that for $q=4$ and 9,
$$
\mbox{``$q$ is a prime pretender just for the bases that are $k$th
roots of 0 or 1 (mod $m$)''}
$$
for a certain $k$  and  $m$.
(It will turn out that such an assertion holds for all the primary
pretenders --- see Table~3.)
They imply that we know the
{\it primary\/} pretender $q_b$ for all but the four residue classes
2, 11, 14, 23 (mod 36):

{\footnotesize
\begin{center}
\begin{tabular}{l@{\,}l@{\hspace{4.5pt}}c@{\hspace{4.5pt}}c@{\hspace{4.5pt}}c
  @{\hspace{4.5pt}}c@{\hspace{4.5pt}}c@{\hspace{4.5pt}}c@{\hspace{4.5pt}}c
  @{\hspace{4.5pt}}c@{\hspace{4.5pt}}c@{\hspace{4.5pt}}c@{\hspace{2pt}}c
  @{\hspace{2pt}}c@{\hspace{2pt}}c@{\hspace{2pt}}c@{\hspace{2pt}}c
  @{\hspace{2pt}}c@{\hspace{2pt}}c@{\hspace{2pt}}c@{\hspace{2pt}}c
  @{\hspace{2pt}}c@{\hspace{2pt}}c@{\hspace{2pt}}c@{\hspace{2pt}}c
  @{\hspace{2pt}}c@{\hspace{2pt}}c@{\hspace{2pt}}c@{\hspace{2pt}}c
  @{\hspace{2pt}}c@{\hspace{2pt}}c@{\hspace{2pt}}c@{\hspace{2pt}}c
  @{\hspace{2pt}}c@{\hspace{2pt}}c@{\hspace{2pt}}c@{\hspace{2pt}}c
  @{\hspace{2pt}}c}
$b$ & $\equiv$ & 0 & 1 & 2 & 3 & 4 & 5 & 6 & 7 & 8 & 9 & 10 & 11 & 12
    & 13 & 14 & 15 & 16 & 17 & 18 & 19 & 20 & 21 & 22 & 23 & 24 & 25
    & 26 & 27 & 28 & 29 & 30 & 31 & 32 & 33 & 34 & 35 \\
$q_b$ & $=$    & 4 & 4 & ? & 6 & 4 & 4 & 6 & 6 & 4 & 4 &  6 &  ? &  4
    &  4 &  ? &  6 &  4 &  4 &  6 &  6 &  4 &  4 &  6 &  ? &  4 &  4
    &  9 &  6 &  4 &  4 &  6 &  6 &  4 &  4 &  6 &  9
\end{tabular}
\end{center}
\normalsize
}

The values of $q_b$ up to 21 for the residue classes mod 1260
missing from the last display are given in Table 1.
In fact $q_b\geq22$ for just the 32 residue classes mod 1260 indicated
by ? in Table~1.
\begin{table}[htb]
\caption{$q_b=10,\,14,\,15,\,21$ for just 108 residue classes mod 1260.}
\begin{center}
\begin{tabular}{r@{\hspace{7pt}}|c@{\hspace{7pt}}c@{\hspace{7pt}}c
  @{\hspace{7pt}}c@{\hspace{7pt}}c@{\hspace{7pt}}c@{\hspace{7pt}}c
  @{\hspace{7pt}}c@{\hspace{7pt}}c@{\hspace{7pt}}c@{\hspace{7pt}}c
  @{\hspace{7pt}}c@{\hspace{6pt}}c@{\hspace{2.4pt}}c@{\hspace{2.4pt}}c
  @{\hspace{2.4pt}}c@{\hspace{2.4pt}}c@{\hspace{2.4pt}}c@{\hspace{2.4pt}}c
  @{\hspace{2.4pt}}c}
$b=$ &  2 & 11 & 14 & 23 & 38 & 47 & 50 & 59 & 74 & 83 & 86 & 95
      & 110 & 119 & 122 & 131 & 146 & 155 & 158 & 167 \\
 \hline \rule{0pt}{12pt}
+0        &  ? & 10 & 14 &  ? &  ? &  ? & 10 & 15 & 15 & 21 & 10 & 10
      &  10 &  14 &  ?  &  10 &  10 &  10 &  ?  &  21 \\
+180      & 14 & 10 & 15 & 14 & 14 &  ? & 10 & 14 & 15 &  ? & 10 & 10
      &  10 &  15 &  14 &  10 &  10 &  10 &  ?  &  ?  \\
+360      &  ? & 10 & 15 &  ? & 21 & 14 & 10 & 15 & 14 &  ? & 10 & 10
      &  10 &  15 &  21 &  10 &  10 &  10 &  14 &  ?  \\
+540      &  ? & 10 & 14 &  ? &  ? & 21 & 10 & 15 & 15 & 14 & 10 & 10
      &  10 &  14 &  ?  &  10 &  10 &  10 &  ?  &  14 \\
+720      & 14 & 10 & 15 & 14 &  ? &  ? & 10 & 15 & 15 &  ? & 10 & 10
      &  10 &  15 &  ?  &  10 &  10 &  10 &  ?  &  ?  \\
+900      & 21 & 10 & 15 & 21 & 14 &  ? & 10 & 14 & 14 &  ? & 10 & 10
      &  10 &  15 &  14 &  10 &  10 &  10 &  14 &  ?  \\
+1080     &  ? & 10 & 15 &  ? &  ? & 14 & 10 & 15 & 15 & 14 & 10 & 10
      &  10 &  15 &  ?  &  10 &  10 &  10 &  21 &  14 \\
\end{tabular}
\end{center}
\end{table}

The number of distinct values of $q_b$ is bounded, since the
 Carmichael number 561 will always
serve if no smaller exponent has been found.  The other numbers which
occur are products of just two prime factors: twice the primes from 2
 to 277; thrice the primes from 3 to
181; five times those primes which are $\equiv1$ (mod 4) from 5 to 109;
seven times those primes which are $\equiv1$ (mod 3) from 7 to 79;
eleven times 11, 31 \& 41; thirteen times 13 \& 37; and the squares of
17, 19 \& 23.

Computer calculations of the numbers in the missing residue classes
for values of $b$ up to 50000 appear in Table 2; the numbers
at the left show the multiples of 1260 to be added.
The programs used to calculate Tables~2 and 3 were straightforward, essentially
using brute force.
\begin{table}[htb]
\caption{$q_b\geq22$ for 32 residue classes mod 1260, $2\leq b\leq51602$.}
\begin{center}
{\scriptsize
\begin{tabular}{r|@{\hspace{0.1pt}}c@{\hspace{0.1pt}}c@{\hspace{0.1pt}}c
  @{\hspace{0.1pt}}c@{\hspace{0.1pt}}c@{\hspace{0.1pt}}c@{\hspace{0.1pt}}c
  @{\hspace{0.1pt}}c@{\hspace{0.1pt}}c@{\hspace{0.1pt}}c@{\hspace{0.1pt}}c
  @{\hspace{0.1pt}}c@{\hspace{0.1pt}}c@{\hspace{0.1pt}}c@{\hspace{0.1pt}}c
  @{\hspace{0.1pt}}c@{\hspace{0.1pt}}c@{\hspace{0.1pt}}c@{\hspace{0.1pt}}c
  @{\hspace{0.1pt}}c@{\hspace{0.1pt}}c@{\hspace{0.1pt}}c@{\hspace{0.1pt}}c
  @{\hspace{0.1pt}}c@{\hspace{0.1pt}}c@{\hspace{0.1pt}}c@{\hspace{0.1pt}}c
  @{\hspace{-0.2pt}}c@{\hspace{-0.6pt}}c@{\hspace{-0.6pt}}c@{\hspace{-0.6pt}}c
  @{\hspace{-0.6pt}}c}
$b=$    & \tiny  2 & \tiny 23 & \tiny 38 & \tiny 47 & \tiny 122 &
 \tiny 158 & \tiny 227 & \tiny 263 & \tiny 338 & \tiny 347 & \tiny 362
 & \tiny 383 & \tiny 443 & \tiny 527 & \tiny 542 & \tiny 563 & \tiny 578
 & \tiny 662 & \tiny 698 & \tiny 758 & \tiny 767 & \tiny 803
 & \tiny 842 & \tiny 878 & \tiny 887 & \tiny 947 & \tiny 983 & \tiny 1067
 & \tiny 1082 & \tiny 1103 & \tiny 1118 & \tiny 1202 \\
 \hline \rule{0pt}{9pt}
0       & 341 & 22 & 38 & 46 &  22 & 158 &  49 &  33 &  26 &  87 &  33
 & 382 &  25 &  33 &  91 &  91 &  34 &  39 &  34 &  33 &  26 &  22
 &  58 & 259 &  91 &  22 &  65 &  22  &  25  &  38  &  25  &  169 \\
1       &  26 & 91 & 22 & 25 &  25 &  25 &  65 & 185 &  34 &  22 &  91
 &  25 &  26 &  38 &  34 &  57 &  22 &  34 &  22 &  25 &  39 &  91
 &  22 &  49 &  38 &  25 &  25 &  26  &  33  &  34  &  39  &  46 \\
2       &  26 & 25 & 91 & 34 &  38 &  26 &  82 &  22 & 113 &  94 &  22
 &  33 &  39 &  22 & 145 &  46 &  38 &  25 &  25 &  22 &  38 &  22
 &  49 &  33 &  25 &  51 &  34 &  22 &   26 &   91 &   34  &  51 \\
3       &  25 & 49 & 22 & 33 &  49 &  22 &  25 &  25 &  25 &  85 &  38
 &  46 &  33 &  25 &  33 &  25 &  91 &  91 &  22 &  26 &  91 &  58
 &  69 &  34 &  26 &  34 &  22 &  74 &   22 &   33 &   62  &  25 \\
4       & 121 &122 & 49 & 85 &  26 &  46 &  46 &  22 &  33 &  65 &  22
 &  22 &  91 &  22 &  25 &  26 &  25 &  49 &  38 &  22 &  25 &  25
 &  25 &  22 &  39 &  82 &  38 &  25 &   39 &   25 &   85  &  49 \\
5       &  33 & 51 &133 & 22 &  26 &  22 &  26 &  34 &  65 &  34 & 133
 &  26 &  22 & 145 &  22 &  33 &  26 &  33 &  91 &  39 &  57 &  65
 &  65 &  74 &  85 & 133 &  22 &  58 &   22 &   22 &   25  &  22 \\
6       &  38 & 34 & 51 & 25 &  25 &  25 &  26 &  91 &  22 &  25 &  34
 &  22 &  65 &  26 &  91 &  62 &  26 & 365 &  46 &  25 &  22 &  51
 &  62 &  22 &  33 &  25 &  25 &  38 &   49 &   57 &   33  &  26 \\
7       &  22 & 25 & 34 & 22 &  33 &  65 &  39 &  38 &  38 &  85 &  25
 &  39 &  22 &  26 &  22 &  22 &  39 &  22 &  25 &  91 &  51 &  34
 &  91 &  26 &  25 &  33 &  26 & 133 &   91 &   22 &   38  &  22 \\
8       &  25 &561 & 25 & 26 &  57 &  58 &  22 &  25 &  22 &  26 &  46
 & 327 &  34 &  25 &  26 &  25 &  33 &  46 &  26 &  91 &  22 & 341
 &  33 &  39 &  22 &  74 &  26 & 142 &   65 &   91 &   22  &  25 \\
9       &  22 & 22 & 62 & 39 &  22 & 121 &  86 &  82 &  51 &  26 & 141
 &  38 &  65 &  34 &  25 &  22 &  25 &  22 &  26 &  46 &  25 &  25
 &  25 &  25 &  85 &  22 &  39 &  25 &  123 &   25 &   91  &  85 \\
10      &  65 & 26 & 33 & 51 &  49 &  49 &  22 &  38 &  34 &  22 &  26
 &  91 &  25 &  91 &  34 &  49 &  22 &  38 &  33 &  38 &  82 &  26
 &  22 &  46 &  22 &  26 &  34 &  51 &   25 &   26 &   22  &  74 \\
11      &  58 & 22 & 26 & 25 &  22 &  25 &  65 &  33 &  62 &  25 &  26
 &  25 &  91 &  33 &  38 & 205 &  65 &  26 &  58 &  25 &  69 &  22
 &  39 &  51 &  65 &  22 &  25 &  22 &   51 &   26 &   34  &  34 \\
12      &  49 & 25 & 22 & 58 & 145 &  33 &  85 &  91 &  26 &  22 &  25
 &  46 &  49 &  65 &  49 &  65 &  22 &  25 &  22 &  34 &  26 &  38
 &  22 &  34 &  25 &  39 &  33 &  57 &   33 &   39 &   26  &  38 \\
13      &  25 & 91 & 25 & 65 &  58 &  46 &  25 &  22 &  25 &  49 &  22
 &  33 &  26 &  22 &  91 &  25 &  62 &  39 &  91 &  22 &  26 &  22
 &  34 &  33 &  49 &  49 &  65 &  22 &   38 &   93 &   26  &  25 \\
14      &  26 & 34 & 22 & 33 &  91 &  22 &  34 &  49 &  39 &  34 &  49
 &  65 &  26 &  62 &  25 &  38 &  25 & 451 &  22 &  85 &  25 &  25
 &  25 &  25 & 118 &  91 &  22 &  25 &   22 &   25 &   39  &  33 \\
15      &  26 & 57 & 34 &169 &  46 &  26 &  62 &  22 &  33 &  38 &  22
 &  22 &  25 &  22 & 145 &  74 &  85 &  62 &  82 &  22 &  33 &  34
 &  38 &  22 &  26 &  91 & 118 &  39 &   25 &   91 &   25  &  38 \\
16      &  33 & 85 & 38 & 22 &  25 &  22 &  38 &  26 &  74 &  25 &  62
 &  25 &  22 &  91 &  22 &  26 &  91 &  33 &  51 &  25 &  34 &  94
 &  49 &  91 &  26 &  25 &  22 &  87 &   22 &   22 &   91  &  22 \\
17      &  62 & 25 & 65 & 85 &  26 &  49 &  33 &  39 &  22 &  65 &  25
 &  22 &  34 &  34 &  65 &  26 &  34 &  25 &  25 &  26 &  22 &  82
 & 123 &  22 &  25 & 106 &  46 &  85 &   39 &  133 &   33  & 133 \\
18      &  22 & 33 & 25 & 22 &  26 &  38 &  25 &  25 &  25 &  57 &  85
 &  26 &  22 &  25 &  22 &  22 &  26 &  22 &  91 &  39 &  38 &  46
 &  38 &  85 & 133 &  33 &  51 &  49 &   85 &   22 &   65  &  22 \\
19      &  49 & 62 & 38 & 34 &  39 &  87 &  22 &  91 &  22 &  33 &  38
 &  26 &  49 &  26 &  25 & 133 &  25 & 146 &  91 & 106 &  22 &  25
 &  25 &  25 &  22 &  46 &  34 &  25 &   57 &   25 &   22  &  26 \\
20      &  22 & 22 & 65 & 26 &  22 &  62 &  39 &  58 &  85 &  49 &  91
 &  39 &  25 &  26 &  85 &  22 &  39 &  22 &  38 &  51 &  46 &  33
 &  58 &  26 &  38 &  22 &  26 &  33 &   25 &   58 &   25  &  34 \\
21      &  94 &142 & 33 & 25 &  25 &  25 &  22 &  49 &  65 &  22 &  49
 &  25 &  57 &  39 &  26 &  91 &  22 &  74 &  26 &  25 &  38 &  49
 &  22 &  39 &  22 &  25 &  25 &  85 &   65 &   86 &   22  & 278 \\
22      &  38 & 22 & 57 & 39 &  22 & 106 &  51 &  33 &  49 &  26 &  25
 &  91 &  65 &  33 &  26 &  65 &  91 &  25 &  25 &  33 &  49 &  22
 &  26 &  49 &  25 &  22 &  39 &  22 &   85 &   38 &   91  &  91 \\
23      &  25 & 26 & 22 & 46 &  65 &  33 &  25 &  25 &  25 &  22 &  26
 &  85 &  85 &  25 &  39 &  25 &  22 & 121 &  22 &  69 &  91 &  26
 &  22 &  91 &  57 &  26 &  33 &  91 &   33 &   26 &   91  &  25 \\
24      &  51 & 39 & 26 & 38 &  34 &  49 &  65 &  22 & 133 &  82 &  22
 &  33 &  46 &  22 &  25 &  49 &  25 &  26 &  49 &  22 &  25 &  22
 &  25 &  25 &  34 &  26 &  91 &  22 &   65 &   25 &  133  &  46 \\
25      &  34 & 58 & 22 & 33 & 145 &  22 &  58 &  46 &  26 &  91 &  39
 &  38 &  25 &  51 &  33 &  34 &  51 &  26 &  22 &  65 &  26 &  62
 & 206 &  65 &  91 &  39 &  22 &  38 &   22 &   33 &   25  &  33 \\
26      &  49 & 74 & 39 & 25 &  25 &  25 & 301 &  22 &  26 &  25 &  22
 &  22 &  26 &  22 &  49 &  91 &  34 &  34 &  57 &  22 &  26 &  91
 &  69 &  22 &  85 &  25 &  25 &  26 &   86 &   51 &   26  &  85 \\
27      &  26 & 25 &145 & 22 &  57 &  22 &  46 & 121 &  34 &  49 &  25
 &  65 &  22 &  46 &  22 &  33 &  57 &  25 &  25 &  85 &  39 &  46
 &  91 & 341 &  25 &  49 &  22 &  26 &   22 &   22 &   39  &  22 \\
28      &  25 & 38 & 25 & 34 &  62 &  26 &  25 &  25 &  22 &  46 &  49
 &  22 &  39 &  25 & 133 &  25 &  46 &  85 & 133 &  74 &  22 &  49
 &  65 &  22 &  26 &  34 & 133 &  34 &   26 &   91 &   33  &  25 \\
29      &  22 & 33 & 91 & 22 &  33 &  39 &  62 &  26 &  49 &  87 &  65
 &  51 &  22 &  86 &  22 &  22 &  25 &  22 &  46 &  26 &  25 &  25
 &  25 &  25 &  26 &  33 &  65 &  25 &   26 &   22 &   38  &  22 \\
30      &  82 &218 & 65 & 49 &  26 & 529 &  22 &  34 &  22 &  33 &  62
 &  34 &  25 &  91 &  38 &  26 &  33 &  65 &  65 &  26 &  22 &  85
 &  33 &  91 &  22 &  85 &  91 &  91 &   25 &  106 &   22  &  85 \\
31      &  22 & 22 &226 & 25 &  22 &  25 &  26 &  87 &  65 &  25 &  46
 &  25 &  91 &  62 &  91 &  22 &  26 &  22 &  49 &  25 &  65 &  33
 &  65 &  38 &  86 &  22 &  25 &  33 &  169 &   86 &   65  &  26 \\
32      &  91 & 25 & 33 & 74 &  39 &  34 &  22 &  57 &  58 &  22 &  25
 &  26 &  65 &  26 &  58 & 217 &  22 &  25 &  25 &  38 & 321 &  34
 &  22 &  26 &  22 &  57 & 206 &  49 &   38 &   46 &   22  &  26 \\
33      &  25 & 22 & 25 & 26 &  22 &  65 &  25 &  25 &  25 & 185 &  33
 &  39 &  49 &  25 &  49 &  25 &  39 &  65 &  51 &  33 &  34 &  22
 &  91 &  26 &  34 &  22 &  26 &  22 &   74 &   85 &   49  &  25 \\
34      &  62 & 57 & 22 & 26 &  91 &  33 & 111 &  46 &  65 &  22 &  91
 & 205 &  34 &  34 &  25 &  91 &  22 &  51 &  22 &  91 &  25 &  25
 &  22 &  25 &  49 &  49 &  26 &  25 &   33 &   25 &   62  &  38 \\
35      &  85 & 26 & 38 & 39 &  91 & 133 &  38 &  22 &  34 &  26 &  22
 &  33 &  25 &  22 &  26 &  65 &  86 &  34 &  26 &  22 &  85 &  22
 &  26 &  33 &  69 & 202 &  39 &  22 &   25 &   34 &   25  &  91 \\
36      &  65 & 26 & 22 & 25 &  25 &  22 &  46 &  86 &  49 &  25 &  26
 &  25 &  33 &  85 &  33 &  38 &  85 &  69 &  22 &  25 &  49 &  26
 &  26 &  49 &  38 &  25 &  22 &  34 &   22 &   26 &   34  &  33 \\
37      &  46 & 25 & 26 & 46 &  86 &  38 &  65 &  22 &  33 &  46 &  22
 &  22 &  38 &  22 &  85 &  58 &  65 &  25 &  25 &  22 &  33 &  39
 &  38 &  22 &  25 &  26 &  91 &  86 &   34 &   26 &  169  &  34 \\
38      &  25 & 49 & 25 & 22 &  46 &  22 &  25 &  25 &  25 &  58 &  38
 &  34 &  22 &  25 &  22 &  25 &  87 &  26 &  49 &  34 &  26 & 178
 &  34 &  65 &  74 &  39 &  22 &  85 &   22 &   22 &   26  &  22 \\
39      & 133 & 51 & 39 & 65 &  65 & 133 &  33 &  34 &  22 &  34 &  87
 &  22 &  26 &  38 &  25 &  46 &  25 &  39 &  38 &  82 &  22 &  25
 &  25 &  22 &  33 & 158 &  38 &  25 &   85 &   25 &   25  &  49 \\
40      &  22 & 33 & 51 & 22 &  33 &  26 &  34 &  58 &  39 &  62 &  34
 &  49 &  22 & 133 &  22 &  22 &  38 &  22 &  58 &  91 &  38 &  51
 &  57 &  94 & 169 &  33 &  46 &  26 &   25 &   22 &   25  &  22
\end{tabular}
\normalsize}
\end{center}
\end{table}

Our final table, Table 3, shows how long it takes before any
particular value of $q_b$ appears; it can be summarized as follows.
The value of $q_b$ is \\
\begin{center}
\begin{tabular}{ccccc}
          &  4       & if $b\equiv$ &    0,1   & (mod 4) \\
     else &  6       & if $b\equiv$ &    0,1   & (mod 3) \\
     else &  9       & if $b\equiv$ &     8    & (mod 9) \\
     else & $\ldots$ & $\ldots$     & $\ldots$ & $\ldots$ \\
     else & 561      & if $b\equiv$ &     0    & (mod 1) \\
\end{tabular}
\end{center}
where the various statements can all be put into the form
$$\mbox{``else $q$ if $b$ is a $k$th root of 0 or 1 (mod $m$)''}$$
for appropriate values of $q$, $k$ and $m$.  The table also gives
the {\bf first base}, that is the least $b$ for which $q_b=q$, and
the {\bf rarity} $r$ of $q$, meaning that $q$ is the primary
pretender for 1 in every $r$ bases.  For example
$$\mbox{25\qquad 4th(25)\qquad 443\qquad 240.62}$$
means that 25 is the primary pretender for the bases that are
4th roots of 0 or 1 $\pmod{25}$ that have not already been coped with,
that the first such base is 443, and that 1 in every 240.62 bases
 has 25 for its primary pretender (in fact 16 in every 3465 bases).

Another example is `else 169 if $b^{12}\equiv0$ or 1 (mod 169)',
i.e., if $b\equiv\pm19^e$ (mod 169),
for $1\leq e\leq6$ where the cases $e=6$ ($b\equiv\pm1$), $e=3$
($b\equiv\pm70$), and $e=2$ or 4 ($b\equiv\pm23$ or $\pm22$) have
already been preempted by $q_b=26$ or 39, by 65, and by 91 respectively.
\clearpage

The largest first base is 10009487, for  $q = 453$, while the greatest rarity
is that of  $q = 519$.

\subsection*{Reference}
The paper was prompted by the table of pseudoprimes to
various bases given by Albert H.~Beiler on p.~42 of his 
{\it Recreations in the Theory of Numbers}, Dover, New York, 1964.

\begin{flushleft}
\begin{tabular}{l}
J. H. Conway \\
Dept.\ of Mathematics, \\
Princeton University, \\
Princeton NJ, 08544 \\
conway@math.princeton.edu \\
~ \\
R. K. Guy \\
Dept.\ of Mathematics, \\
The University of Calgary, \\
Calgary, Alberta, \\
Canada T2N 1N4. \\
rkg@cpsc.ucalgary.ca \\
~ \\
W. A. Schneeberger \\
Dept.\ of Mathematics, \\
Princeton University, \\
Princeton NJ, 08544 \\
william@math.princeton.edu \\
~\\
N. J. A. Sloane \\
AT\&T Bell Laboratories \\
Mathematical Sciences Research Center \\
Murray Hill NJ, 07974 \\
njas@research.att.com
\end{tabular}
\end{flushleft}

\begin{table}[htb]
\caption{The first base and rarity of the 132 primary pretenders.}
\begin{center}
{\footnotesize
\begin{tabular}{r@{\hspace{4pt}}l@{\hspace{-2pt}}r@{\hspace{2pt}}
   r|r@{\hspace{2pt}}l@{\hspace{0pt}}r@{\hspace{2pt}}
   r|r@{\hspace{2pt}}l@{\hspace{-2pt}}r@{\hspace{2pt}}r}
$q$ & roots      & first  & rarity & $q$\verb+ + & roots      & first
 & rarity & $q$\verb+ + & roots      & first  & rarity  \\
    & $k$th($m$) & base   & one in &             & $k$th($m$) & base
 & one in &             & $k$th($m$) & base   & one in  \\
 \hline  \rule{0pt}{9pt}
  4&1st(4)   &       0&	    2\ \ \ \ &	159&2nd(53)  &   94763& 
 83341.92 &	361&18th(361)&   58727&   159201.47  \\
  6&1st(3)   &       3&	    3\ \ \ \ &	166&1st(83)  &  247838& 
   69173.80 &	362&1st(181) & 1050887&   800429.64  \\
  9&2nd(9)   &      26&	   18\ \ \ \ &	169&12th(169)&    1202& 
   22203.93 &	365&4th(73)  &    8222&   313017.17  \\
 10&1st(5)   &      11&	  22.5\ \ &	177&2nd(59)  &  111863& 
  105468.69 &	381&2nd(127) &  923162&  1150798.45  \\
 14&1st(7)   &      14&	  52.5\ \ &	178&1st(89)  &   48683& 
   83809.94 &	382&1st(191) &     383&   886300.41  \\
 15&2nd(5)   &      59&	   63\ \ \ \ &	183&2nd(61)  &  186842& 
  113673.26 &	386&1st(193) &  470342&   905058.10  \\
 21&2nd(7)   &      83&	 157.5\ \ &	185&4th(37)  &    1523& 
   33318.02 &	393&2nd(131) &  480638&  1222539.21  \\
 22&1st(11)  &      23&	     216.56 &	194&1st(97)  &   58298& 
  100995.26 &	394&1st(197) &  384347&   940782.12  \\
 25&4th(25)  &     443&	     240.62 &	201&2nd(67)  &   86027& 
  138204.04 &	398&1st(199) &  278402&   960080.22  \\
 26&1st(13)  &     338&	     391.01 &	202&1st(101) &   45047& 
  109051.62 &	411&2nd(137) &  786242&  1315845.99  \\
 33&2nd(11)  &     263&	     639.84 &	205&4th(41)  &   14423& 
   41858.20 &	417&2nd(139) &  303158&  1345305.23  \\
 34&1st(17)  &     578&	     679.83 &	206&1st(103) &   32342& 
  119760.96 &	422&1st(211) &  231467&  1043600.74  \\
 38&1st(19)  &      38&	     861.12 &	213&2nd(71)  &   53462& 
  163633.79 &	427&6th(61)  &  149558&   139812.54  \\
 39&2nd(13)  &     662&	    1114.39 &	214&1st(107) &   79502& 
  128741.29 &	445&4th(89)  &  739022&   462456.86  \\
 46&1st(23)  &      47&	    1281.55 &	217&6th(31)  &   40883& 
   17165.50 &	446&1st(223) &  592958&  1227712.86  \\
 49&6th(49)  &     227&	     854.36 &	218&1st(109) &   37823& 
  155920.00 &	447&2nd(149) &  141698&  1633246.98  \\
 51&2nd(17)  &    3467&	    2135.92 &	219&2nd(73)  &  169067& 
  206921.87 &	451&10th(41) &   18302&    50339.80  \\
 57&2nd(19)  &    1823&	    2593.62 &	226&1st(113) &   39098& 
  167015.51 &	453&2nd(151) &10009487&  2143037.38  \\
 58&1st(29)  &     842&	    2350.46 &	237&2nd(79)  &  141962& 
  231715.22 &	454&1st(227) &  283523&  1643477.99  \\
 62&1st(31)  &    4898&	    2698.68 &	249&2nd(83)  &  357563& 
  246959.64 &	458&1st(229) &  277778&  1672695.37  \\
 65&4th(13)  &     983&	     930.58 &	254&1st(127) &  232538& 
  196024.21 &	466&1st(233) &  860702&  1716907.59  \\
 69&2nd(23)  &    4622&	    4885.55 &	259&6th(37)  &     878& 
   25091.09 &	469&6th(67)  &  473987&   237840.01  \\
 74&1st(37)  &    4847&	    4519.13 &	262&1st(131) &  162047& 
  234781.00 &	471&2nd(157) & 2264567&  2457680.13  \\
 82&1st(41)  &    2747&	    5293.84 &	265&4th(53)  &   98663& 
   91000.38 &	478&1st(239) & 6085658&  1907095.94  \\
 85&4th(17)  &    4127&	    1900.35 &	267&2nd(89)  &  232823& 
  329876.41 &	481&12th(37) &  108803&   112655.45  \\
 86&1st(43)  &   11567&	    6809.60 &	274&1st(137) &  112478& 
  262750.39 &	482&1st(241) & 2252387&  2262497.08  \\
 87&2nd(29)  &     347&	    8968.74 &	278&1st(139) &   27662& 
  270535.59 &	485&4th(97)  &  968567&   889852.41  \\
 91&6th(13)  &     542&	     689.90 &	289&16th(289)&  197138& 
   67140.22 &	489&2nd(163) & 3166763&  3114483.44  \\
 93&2nd(31)  &   17483&    20007.20 &	291&2nd(97)  &  124547& 
  398645.06 &	501&2nd(167) & 4881242&  3211811.04  \\
 94&1st(47)  &    2867&	   16791.76 &	298&1st(149) &  142742& 
  315947.41 &	502&1st(251) & 1738427&  2457818.82  \\
106&1st(53)  &   22367&	   19776.96 &	301&6th(43)  &   32987& 
   42986.04 &	505&4th(101) & 2128262&   967334.31  \\
111&2nd(37)  &   43067&	   27144.85 &	302&1st(151) &  150698& 
  360605.13 &	511&6th(73)  &  210962&   342597.56  \\
118&1st(59)  &   18527&	   23552.15 &	303&2nd(101) &  485102& 
  479193.40 &	514&1st(257) & 2338187&  2751486.73  \\
121&10th(121)&    5042&     9090.30 &	305&4th(61)  &  287138& 
  141802.12 &	519&2nd(173) & 1150103&  3690229.26  \\
122&1st(61)  &    5063&    27725.42 &	309&2nd(103) &   79103& 
  511500.54 &	526&1st(263) &  256163&  2854500.87  \\
123&2nd(41)  &   12422&    36653.95 &	314&1st(157) &  115238& 
  401527.92 &	529&22nd(529)&   37958&   503092.10  \\
129&2nd(43)  &   66047&    39547.68 &	321&2nd(107) &   41087& 
  544005.57 &	537&2nd(179) & 7345622&  4047604.68  \\
133&6th(19)  &    2858&     3954.76 &	326&1st(163) &  296987& 
  426312.06 &	538&1st(269) & 2735462&  3093197.89  \\
134&1st(67)  &   87302&    44161.58 &	327&2nd(109) &   10463& 
  566650.81 &	542&1st(271) &  183467&  3139537.94  \\
141&2nd(47)  &   11702&    61146.80 &	334&1st(167) &  127922& 
  446371.15 &	543&2nd(181) & 4503098&  4178269.82  \\
142&1st(71)  &   11147&    49334.35 &	339&2nd(113) &  851567& 
  600572.10 &	545&4th(109) & 4453598&  1244091.57  \\
145&4th(29)  &    3062&    18589.75 &	341&10th(31) &      2 & 
   16379.23 &	553&6th(79)  &  281738&   454571.92  \\
146&1st(73)  &   24602&    56543.84 &	346&1st(173) &  846662& 
  708402.09 &	554&1st(277) &  581423&  3497678.40  \\
158&1st(79)  &     158&    62914.98 &	358&1st(179) &  257402& 
  741543.71 &	561&1st(1)   &   10103&    25437.66 
\end{tabular}
\normalsize}
\end{center}
\end{table}
\end{document}